\documentclass[10pt]{article}
\usepackage{amsmath}
\usepackage{amsthm}
\usepackage{amssymb}
\usepackage{amscd}
\usepackage{amsfonts}
\usepackage{enumerate}
\usepackage{amsbsy}
\usepackage[all]{xy}
\usepackage{graphicx}
\usepackage{shapepar}
\usepackage[latin1]{inputenc}

\def\H{{\mathbb H}}
\def\L{{\mathbb L}}

\def\F{{\mathbb F}}

\def\K{{\mathbb K}}
\def\IK{{\mathfrak I}_{\K}}
\def\IF{{\mathfrak I}_{\F}}
\def\Q{{\mathbb Q}}
\newcommand{\Oo}[1]{{\cal O}_{#1}}

\def\Z{{\mathbb Z}}
\def\N{{\mathbb N}}
\def\q{{\mathfrak q}}
\def\p{{\mathfrak p}}

\def\OF{\Oo{\F}}

\def\GF{{Cl_{\F}}}

\def\G2{{Cl_{2}}}

\def\OK{\Oo{\K}}

\newcommand{\ord}[2]{{\rm ord}_{#1}\left(#2\right)}

\newcommand{\leg}[2]{{\dps\left(\dfrac{#1}{#2}\right)}}

\newcommand{\ide}[1]{\langle{#1}\rangle}
\newcommand{\ideF}[1]{\left\langle{#1}\right\rangle_{\F}}

\newcommand{\ideK}[1]{\langle{#1}\rangle_{\K}}
\newcommand{\mcd}[1]{{\rm g.c.d.}\left(#1\right)}

\newcommand{\nkf}[1]{{N_{\K/\F}\left(#1\right)}}
\newcommand{\nk}[1]{{N_{\K/\Q}\left(#1\right)}}
\newcommand{\nf}[1]{{N_{\F/\Q}\left(#1\right)}}

\def\JK{{{{\mathfrak J}}_\K}}
\def\Irr{{\rm Irr}}

\def\JF{{{{\mathfrak{J}}_\F}}}
\def\C{{\mathbb C}}

\newcommand{\ideL}[1]{\left\langle{#1}\right\rangle_{\L}}

\def\JJ{{{{\mathfrak{J}}}}}
\def\II{{{{\mathfrak{I}}}}}

\newcounter{conta}
\newtheorem {teo} {Theorem}
\newtheorem {prop} [teo] {Proposition}
\newtheorem {cor} [teo] {Corollary}
\newtheorem {lem} [teo] {Lemma}

\def\dps{\displaystyle}

\def\gal{\mbox{\rm Gal}}

\def\GK{{Cl_{\K}}}

\title{ An infinite family of pure quartic fields with class number $\equiv 2\pmod{4}$}

\author{Alejandro Aguilar-Zavoznik\\
\small{aaz@correo.azc.uam.mx}\\
\small Departamento de Ciencias B\'asicas\\
\small Universidad Aut\'onoma Metropolitana-Azcapotzalco\\
\small Av. San Pablo No. 180, Col. Reynosa Tamaulipas\\
\small C.P. 02200 Del. Azcapotzalco M\'exico D.F.\\
Mario Pineda-Ruelas\\
\small{mpr@xanum.uam.mx}\\
\small{Departamento de Matem\'aticas}\\
\small{Universidad Aut\'onoma Metropolitana-Iztapalapa}\\
\small{Av. San Rafael Atlixco No. 186, Col. Vicentina}\\
\small{C.P. 09340 Del. Iztapalapa M\'exico D.F. }\\
}

\begin{document}
\maketitle

\begin{abstract}

Let us consider the pure quartic fields of the form $\K=\Q(\sqrt[4]{p})$ where 
$0<p\equiv 7\pmod{16}$ is a prime integer.
We prove that the $2$-class group of $\K$ has order $2$. As a consequence of this, if the class number of $\K$ is $2$, 
then the Hilbert class field of $\K$ is $\H_\K=\K(\sqrt{2})$. Finally, we find a criterion to decide if an ideal of the ring of integers
or $\K$ is principal or non-principal.

 \let\thefootnote\relax
\end{abstract}

\section{Introduction}

The main goal of this paper is to prove:
\begin{teo}\label{mainresult}
 Let $\K=\Q(\sqrt[4]{p})$ with $0<p\equiv 7\pmod{16}$ a prime integer. 
 Then, the class number $h_\K\equiv 2\pmod{4}$. Equivalently, the 
 $2$-class group of $\K$ is isomorphic to $\Z/2\Z$.
\end{teo}

The $2$-class group of number fields has been widely studied, Gauss described the $2$-rank of the class
group of a quadratic number field using the language of binary quadratic forms. Since then, a lot of work 
has been done on this subject, if the reader is interested in this subject, we suggest the next bibliography: 
\cite{azpr3}, \cite{azpr4}, \cite{6}, \cite{4}, \cite{3}, %8
 \cite{basilla}, \cite{8}, \cite{a7}, \cite{bosma}, \cite{1}, \cite{a1},   \cite{hasse}, 
 \cite{a6}, \cite{7},  \cite{a3},  \cite{5}, \cite{mollin1}, \cite{mollin2}, \cite{a9},  \cite{nakano}, \cite{2}, %26     
\cite{shanks}. 

An important result is the Ambiguous Class Number Formula that 
states that, if $\K/\F$ is a cyclic extension, then 
$${\rm Am}(\K/\F)=h_\F\dfrac{\prod e(\p)}{(\K:\F) (E_\F:E_\F\cap \nkf{\K^\times})}$$ 
where the product is over all primes (finite and infinite) and $E_\F$ is the group of units
of  the ring of integers of $\F$ (see \cite{lemmermeyer} Theorem 11.14).
If $\K/\F$ is a quadratic extension, ${\rm Am}(\K/\F)$ helps us to find the $2$-rank of the 
ideal class group of $\K$. In particular, in \cite{8}, \cite{a7}  the authors study
families of number fields with given $2$-class group. The purpose of this work is to show
that the $2$-class group of a family of pure quartic fields has order $2$. 
% a family
% study the previous problem 
% in a family of pure quartic fields.

Let $p$ and $\K$ be as in Theorem \ref{mainresult}, $\F=\Q(\sqrt{p})$ and $\H_\K$ the Hilbert class field
of $\K$. We will denote $\OF$, $h_\F$, $\GF$, $\G2$ and ${\cal U}_\F$ the ring of integers, class number, class group, $2$-class group and 
group of units of an arbitrary number field $\F$ respectively and if $\F$ is a quadratic number field
$U_\F$ is the fundamental unit of $\F$. Given an extension of number fields $\K/\F$,  let $\nkf{\alpha}$
be the relative norm respect to $\K/\F$ of $\alpha\in\K$, $\nk{\alpha}$ the absolute norm of $\alpha\in\K$,
$\nf{A}$ the absolute norm of $A\in\F$. We will use the symbol $\ideF{A_1,\ldots,A_n}$ to denote the ideal of 
$\OF$ generated by $A_1,\ldots,A_n\in\OF$ and $\ideK{\alpha_1,\ldots,\alpha_n}$ the ideal of $\OK$ generated by
$\alpha_1,\ldots,\alpha_n\in\OK$. The symbol $\ord{a}{b}=c$  means that $a^c$ is the greatest power
of $a$ that divides $b$ and $a^c||b$ denotes that $\ord{a}{b}=c$. We will denote $\N_0=\{0,1,2,\ldots\}$ the set of the natural
numbers starting from $0$.

To prove the Theorem \ref{mainresult} we will use the next facts about the Hilbert Class Field of $\K$.
We know that $\GK\cong C_{n_1}\times\cdots\times C_{n_k}$ where $C_n=\Z/n\Z$ and $n_i=p_i^{e_i}$
for some rational prime $p_i$. If the $2$-rank of $\GK$ is $r$, then there is a group $H\subseteq \GK$ such
that $H\cong C_2^r$ and $r$ is the maximum integer that satisfies this. Since $\GK$ is isomorphic to the 
Galois group of $\H_\K/\K$, then there is another group $H_1\subseteq \gal(\H_\K/\K)$ such that 
$H_1\cong \GK/H$. The Galois Group of the fixed field $\H_\K^{H_1}$ is isomorphic to 
$C_2^r$. So, $r$ is the maximum integer such that there is a non-ramified Galois extension over $\K$ with 
Galois group isomorphic to $C_2^r$.

\section{Previous results}

 In this section we will state some results found in \cite{azpr} and \cite{azpr2} that we will use in the next sections.
 
\begin{teo}\label{resultofazpr2} 
   Let $\K=\Q(\sqrt[4]{p})$ with a rational prime $0<p\equiv 7\pmod{16}$, $\L=\K(\sqrt{\alpha})$ with
  $\alpha=a_1+a_2\sqrt[4]{p}+a_3\sqrt{p}+a_4\sqrt[4]{p^3}\in\OK$ and $\L\neq \K$. Then, $2$ does not ramifies completely in 
  $\L$ if and only if one of the following assertions holds:
  \begin{enumerate}
  \item[\rm 1.] $\alpha=U_\F$ the fundamental unit of $\F=\Q(\sqrt{p})$.
  \item[\rm 2.] $\nk{\alpha}\equiv 4\pmod{8}$, $a_1\equiv a_3\pmod 8$ are odd, $a_2\equiv 2\pmod 4$ and $a_4\equiv 0 \pmod 4$.
  \item[\rm 3.] $\nk{\alpha}\equiv 4\pmod{8}$, $a_1\not\equiv a_3\pmod 4$ odd, $a_1+a_3\equiv 0\pmod 8$, $a_2\equiv 0\pmod{4}$ 
   and $a_4\equiv 2\pmod{4}$.
  \item[\rm 4.] $N(\alpha)$ odd, $\alpha\not\in {\cal U}_\K$, $a_1$ odd, $a_2\equiv a_4\pmod{4}$ even and $a_3\equiv 0\pmod{4}$.
  If  $a_3$ is odd and $a_1,a_2,a_4$ are even we multiply $\alpha$ by $\sqrt{p}$ to get a new $\alpha$ that generates the same field
  $\L$ with $a_1$ odd and $a_2,a_3,a_4$
  even.
%   (see the comment in the subsection ``$\nk{\alpha}$ odd'' before {\rm Proposition \ref{principal4}}).
  \end{enumerate}
%   In any of this cases, there is a $\gamma\in\OK$ such that if $\nk{\gamma}$ is odd then $\ideK{2}=\p^4$ and if 
%   $\nk{\gamma}$ is even then $\ideK{2}=\p_1^4\p_2^4$, where $\p,\p_1,\p_2$ are prime ideals.

%  
%  Let $\K=\Q(\sqrt[4]{p})$ with $0<p\equiv 7\pmod{16}$ a rational prime and
%  $\L=\K(\sqrt{\alpha})$ with $\alpha=a_1+a_2\sqrt[4]{p}+a_3\sqrt{p}+a_4\sqrt[4]{p^3}\in\OK$ and $\L\neq \K$. Then $2$ does not ramifies in $\L/\K$
%  if and only if one of the next conditions holds:
% \begin{enumerate}[\rm 1.]
%  \item $\nk{\alpha}\equiv 4\pmod 8$ with $a_1\equiv a_3\pmod 8$ odd, $a_2\equiv 2\pmod{4}$ and $a_4\equiv 0\pmod{4}$.
%  \item $\nk{\alpha}\equiv 4\pmod 8$ such that $a_1\not\equiv a_3\pmod 4$, $a_1+a_3\equiv 0\pmod 8$, $a_2\equiv 0\pmod 4$ and $a_4\equiv 2\pmod 4$.
%  \item $\nk{\alpha}$ odd, where $a_1$ odd, $a_2\equiv a_4\pmod 4$ even and $a_3\equiv 0\pmod 4$.
%  
% \end{enumerate}
\end{teo}
\begin{dem}
 See \cite{azpr2}, Theorem 1.1.\qed
\end{dem}
  $$\xymatrix{  
     \L=\K(\sqrt{\alpha}) \ar@{-}[d]  \\
   \K=\Q(\sqrt[4]{p})  \\
     \F =\Q(\sqrt{p}) \ar@{-}[u]  \\
  }
  $$

\begin{prop}\label{3prop11}
 Let $0<p\equiv 3\pmod{4}$ be a rational prime number and $\K=\Q(\sqrt[4]{p})$. The prime $2$ ramifies completely in $\K$.
\end{prop}
\begin{dem}
 See \cite{azpr}, Proposition 11. \qed\bigskip
\end{dem}

\begin{prop}\label{3prop12}
 Let $d\equiv 7\pmod{16}$ be a square-free integer, $\K=\Q(\sqrt[4]{d})$ and the ideal $\IK=\ideK{2,1+\sqrt[4]{d}}$. 
 Then, $\IK$ is a non principal ideal.
\end{prop}
\begin{dem}
 See \cite{azpr}, Proposition 12. \qed \bigskip
\end{dem}

We have that the ideal $\IK$ of the previous proposition is non principal, but $\IK^2$ is principal generated by $L_2$, where 
$L_2$ is such that $2=L_2^2U_\F$.

 \begin{prop} \label{3prop13}
  Let $d\in\Z$ be an odd square-free number, $\K=\Q(\sqrt[4]{d})$, $\F=\Q(\sqrt{d})$,
  $f(x)=x^2+A_1x+A_0\in\OF[x]$, $\alpha\in\C$ with $f(\alpha)=0$ and
  $\Delta_f=A_1^2-4A_0$. Then $\alpha\in\OK$ if and only if there is $C\in\OF$ such that $\Delta_f=C^2$ or
  $\Delta_f=C^2\,\sqrt{d}$. In the first case $\alpha\in\OF$, in the second one $\alpha\in\OK-\OF$.
 \end{prop}
 \begin{dem}
  See \cite{azpr}, Proposition 13.\qed \bigskip
 \end{dem}
 
 \begin{prop}\label{3prop14}
 Let $\alpha>0$ be an element of $\OK-\OF$ such that $\nk{\alpha}=B^2$ for some $B\in\OF$ and 
 $f(x)=\Irr(\alpha,\OF)=x^2-A x+B^2$. There exists $C\in\OF$ such that
 $C^2=A\pm2 B$ for one of the signs if and only if $\sqrt{\alpha}\in\OK$.
 \end{prop}
 \begin{dem}
  See \cite{azpr}, Proposition 14.\qed \bigskip
 \end{dem}

 \begin{teo}\label{3teo18}
  Let $0<p\equiv 7\pmod {16}$ be a rational prime number, $\F=\Q(\sqrt{p})$, $U_\F$ be the fundamental
  unit of $\F$, and $\K=\Q(\sqrt[4]{p})$. Then, the group of units of $\OK$ has the form 
  $\ide{-1,\mu_1,\mu_2}$, where $\nkf{\mu_1}=1$ and $|\nkf{\mu_2}|=U_\F$.
 \end{teo}
 \begin{dem}
  See Theorem 18, \cite{azpr}.\qed\bigskip
 \end{dem}

\section{Proof of Theorem 1}

We will find the only non-ramified quadratic extension of $\K$ using Theorem \ref{resultofazpr2}. First we will
see what happens when no ideal $\IK\subseteq\OK$ satisfies $\IK^2=\ideK{\alpha}$.

% Primero vamos a estudiar lo que sucede cuando no existe ning\'un ideal $\IK\subseteq \OK$ tal que $\IK^2=\ideK{\alpha}$.
% 
\begin{prop}
 Let $\K=\Q(\sqrt[4]{p})$ with $p$ a positive rational prime and $\alpha\in\OK$ a square-free element in 
 any of it's factorizations. If $\p\subseteq \OK$ is a prime ideal such that  
 $\ord{\p}{\ideK{\alpha}}$ is odd, then $\p$ ramifies completely in $\L/\K$.
% 
%  Sea $\K=\Q(\sqrt[4]{p})$ con $p$ un primo racional positivo y $\alpha\in\OK$ libre de cuadrados en todas sus factorizaciones.
%  Supongamos que existe un ideal $\p\subseteq \OK$ tal que $\ord{\p}{\ideK{\alpha}}$ es impar. Entonces $\p$ se ramifica totalmente en $\L/\K$.
\end{prop}
\begin{dem}
Consider $\q_\L=\ideL{\p,\sqrt{\alpha}}$. Since $\ord{\p}{\ideK{\alpha}}$ is odd, then there is $t\in\N$ odd such that
$\ideK{\p}^t\mid \mid \ideK{\alpha}$. The ideal $\ideL{\alpha}$ is a square since $\sqrt{\alpha}\in\L$, 
which implies that any prime ideal that divides $\ideL\alpha$ must appear an even number of times
in the factorization of $\alpha$. Since $t$ is odd, $\q_\L^{2t}\mid\mid\ideL{\alpha}$ and 
$\ideL{\p}=\q_\L^2$. Therefore, $\p$ ramifies completely in $\L/\K$.
\qed\bigskip
% 
% 
% Sea $\ideL{\p,\sqrt{\alpha}}$. Como $\ord{\p}{\ideK{\alpha}}$, es impar, entonces existe $t\in\N$ impar tal que
% $\ideK{\p}^t\mid \mid \ideK{\alpha}$. El ideal $\ideL{\alpha}$ es un cuadrado porque $\sqrt{\alpha}\in\L$, lo que implica que 
% cada ideal primo que divide a $\ideL\alpha$ debe estar elevada a una potencia par. Como $t$ es impar, 
% para que $\q_\L^{2t}\mid\mid\ideL{\alpha}$, necesariamente
% $\ideL{\p}=\q_\L^2$. Entonces $\p$ se ramifica totalmente en $\L/\K$.
% \qed
%\bigskip
\end{dem}

As a consequence of the previous result, if $\ideK{\alpha}$ is not a square, then $\L/\K$ is a ramified extension.
So, it remains to see what happens when $\ideK{\alpha}=\IK^2$ for some ideal $\IK\subseteq \OK$.
Remember that using Gauss Theorem on the $2$-rank of the class group of a quadratic field (\cite{mollant}, 
Theorem 3.70), if $\F=\Q(\sqrt{p})$ with a rational prime $p$, then $h_\F$ is odd.

% La proposici\'on anterior nos muestra que si $\ideK{\alpha}$ no es un cuadrado, entonces la extensi\'on $\L/\K$ es ramificada.
% De esta forma nos falta considerar los casos en que $\ideK{\alpha}=\IK^2$ para alg\'un ideal $\IK\subseteq \OK$.
% Notemos que, por el Teorema del 2-rango de Gauss, si $\F=\Q(\sqrt{p})$ con $p$ un primo racional, entonces $h_\F$ es impar.
% 

\begin{lem}
 Let $\K=\Q(\sqrt[4]{p})$ and $\F=\Q(\sqrt{p})$ with $0<p\equiv 7\pmod{16}$ a rational prime. 
 If $\alpha\in\K$ is such that $\ideK{\alpha}=\IK^2$ for some ideal $\IK\subseteq \OK$,
 then there is an element $\beta\in\OK$ such that $\ideK{\beta}=\ideK{\alpha}$
 and $\nkf{\beta}=B^2$ for some $B\in\OF$.
% 
%  Sean $\K=\Q(\sqrt[4]{p})$ y $\F=\Q(\sqrt{p})$ con $p\equiv 7\pmod{16}$ un primo racional positivo. Si $\alpha\in\K$ es tal que 
%  $\ideK{\alpha}=\IK^2$ para alg\'un ideal $\IK\subseteq \OK$, entonces existe $\beta\in\OK$ tal que $\ideK{\beta}=\ideK{\alpha}$
%  y $\nkf{\beta}=B^2$ para alg\'un $B\in\OF$.
\end{lem}
\begin{dem}
 Since %the norm of an ideal is multiplicative and 
 $\ideK{\alpha}=\IK^2$, then $\nkf{\ideK{\alpha}}=\IF^2$,
 for some $\IF\subseteq\OF$.
 Also, there is $B\in \OF$, $\IF=\ideF{B}$ since $h_\F$ is odd.  
 Hence $\nkf{\ideK{\alpha}}=\ideF{\nkf{\alpha}}=\ideF{B^2}$. 
 From the previous equality, $\nkf{\alpha}=B^2\, U$ with $U\in{\cal U}_\F$. We can suppose that  $U=\pm 1$ or $U=\pm U_\F$. 
 If $U=\pm U_\F$, then $\nkf{\alpha/\mu_2}=\pm B^2$, where $\mu_2$ is the generator of ${\cal U}_\K$ with $\nkf{\mu_2}=U_\F$
 (see Theorem \ref{3teo18}).
 Let $\beta=\alpha$ or $\beta=\alpha/\mu_2$ such that $\nkf{\beta}=\pm B^2$. 

 In  $\OF$, $(a_1+a_2\sqrt{p})^2=a_1^2+p\,a_2^2+2\,a_1\,a_2\,\sqrt{p}$, hence, 
 the squares modulo $\ideF{\sqrt{p}}$ are the same squares of $\Z$ modulo $p$.
 If $a\in\Z$, then $\leg{a}{p}=1$ implies $\leg{-a}{p}=-1$, so, if $A\in\OF$ is a square modulo
 $\ideF{\sqrt{p}}$ then $-A$ is not a square modulo $\ideF{\sqrt{p}}$.
 On the other hand, if $\beta=b_1+b_2\sqrt[4]{p}+b_3\sqrt{p}+b_4\sqrt[4]{p^3}$, we have:
 $$
  \begin{array}{rcl}
   \nkf{\beta} & = & (b_1+b_3\sqrt{p})^2-\sqrt{p}(b_2+b_4\sqrt{p})^2\\
      & = & (b_1^2+b_3^2\,p-2\,p\,b_2\,b_4)+\sqrt{p}(2\,b_1\,b_3-b_2^2-p\,b_4^2).\\
  \end{array}
 $$
 This shows that $\nkf{\beta} \equiv b_1^2 \pmod{ \ideF{\sqrt{p}}}$ is a square modulo $\ideF{\sqrt{p}}$, therefore
 $\nkf{\beta}=B^2$.\qed \bigskip 
\end{dem}

% \begin{dem}
%  Como la norma de ideales es multiplicativa y $\ideK{\alpha}=\IK^2$, entonces $\nkf{\ideK{\alpha}}=\IF^2$,
%  para alg\'un $\IF\subseteq\OF$.
%  Adem\'as, para cierto $B\in \OF$, $\IF=\ideF{B}$ ya que $h_\F$ es impar.  
%  As\'i $\nkf{\ideK{\alpha}}=\ideF{\nkf{\alpha}}=\ideF{B^2}$. 
%  De la igualdad anterior, $\nkf{\alpha}=B^2\, U$ con $U\in{\cal U}_\F$. Podemos suponer $U=\pm 1$ \'o $U=\pm U_\F$. 
%  Si $U=\pm U_\F$, entonces $\nkf{\alpha/\mu_2}=\pm B^2$, donde $\mu_2$ es el generador de ${\cal U}_\K$ con $\nkf{\mu_2}=U_\F$.
%  Sea $\beta=\alpha$ \'o $\beta=\alpha/\mu_2$ tal que $\nkf{\beta}=\pm B^2$. 
% 
%  En  $\OF$, $(a_1+a_2\sqrt{p})^2=a_1^2+p\,a_2^2+2\,a_1\,a_2\,\sqrt{p}$, as\'i, 
%  los cuadrados m\'odulo $\ideF{\sqrt{p}}$ son los mismos cuadrados de $\Z$ m\'odulo $p$.
%  Si $a\in\Z$, entonces $\leg{a}{p}=1$ implica $\leg{-a}{p}=-1$, as\'i que, si $A\in\OF$ es un cuadrado m\'odulo 
%  $\ideF{\sqrt{p}}$ entonces $-A$ no es un cuadrado m\'odulo $\ideF{\sqrt{p}}$.
%  Por otra parte, si $\beta=b_1+b_2\sqrt[4]{p}+b_3\sqrt{p}+b_4\sqrt[4]{p^3}$, tenemos:
%  $$
%   \begin{array}{rcl}
%    \nkf{\beta} & = & (b_1+b_3\sqrt{p})^2-\sqrt{p}(b_2+b_4\sqrt{p})^2\\
%       & = & (b_1^2+b_3^2\,p-2\,p\,b_2\,a_4)+\sqrt{p}(2\,b_1\,b_3-b_2^2-p\,b_4^2).\\
%   \end{array}
%  $$
%  Lo anterior nos muestra que $\nkf{\beta} \equiv b_1^2 \pmod{ \ideF{\sqrt{p}}}$ es un cuadrado m\'odulo $\ideF{\sqrt{p}}$, por lo que
%  $\nkf{\beta}=B^2$.\qed \bigskip 
% \end{dem}

% aquí 5 de sep de 2013.

In the next theorem, assertion 7 is the important one, 1-6 are stated to give a better understanding of the proof.

% La parte importante de la siguiente proposici\'on es la afirmaci\'on 7, los primeros seis incisos sirven como gu\'ia 
% para justificar que $\IK$ es un ideal principal.
% 
\begin{prop}\label{casos1y2}
 Let $\K=\Q(\sqrt[4]{p})$, $\F=\Q(\sqrt{p})$ with $0<p\equiv 7\pmod{16}$ a rational prime,
 $A_1=a_1+a_3\sqrt{p}$, $A_2=a_2+a_4\sqrt{p}$
 and $\alpha=A_1+A_2\sqrt[4]{p}\in\OK-\OF$ such that
 $\ideK{\alpha}=\IK^2$, $\nkf{\alpha}=B^2$ and $\alpha$ satisfies one of the conditions {\rm 1}, 
 {\rm 2} or {\rm 3} of {\rm Theorem \ref{resultofazpr2}}.
 Then the next assertions hold:
 \begin{enumerate}[\rm 1.]
  \item If $B=b_1+b_2\sqrt{p}$, then $b_1,b_2$ have the same parity if 
  {\rm 1} or {\rm 2} of Theorem \ref{resultofazpr2} holds, 
  and if assertion {\rm 3} holds, $b_1$ is odd and $b_2\equiv 0\pmod{4}$.
  \item $\ideF{L_2}^2\mid \mid \ideF{A_1+B}+\ideF{A_1-B}$ where $L_2$ is as defined after Propoposition \ref{3prop12}.
  \item If {\rm 1} or {\rm 2} of Theorem \ref{resultofazpr2} holds, $\ideF{A_1+B}+\ideF{A_1-B}=\ideF{L_2}(\ideF{A_1}+\ideF{B})$. 
  If assertion {\rm 3} holds,
    $\ideF{A_1+B}+\ideF{A_1-B}=\ideF{2}(\ideF{A_1}+\ideF{B})$.
  \item Let $\p_\F$  be a prime ideal of $\OF$ such that $\p_\F\mid \ideF{A_1}$ and $\p_\F\mid \ideF{B}$. If $\p_\F$ is inert in $\K/\F$ and 
  $\p_\F^k\mid\mid \ideF{A_1}+\ideF{B}$, then $k$ is even.
  \item Let $\p_\F$  be a prime ideal of $\OF$ such that $\p_\F\mid \ideF{A_1}$ and $\p_\F\mid \ideF{B}$. If $\p_\F$ splits in
   $\K/\F$ and $\p_\F^k\mid\mid \ideF{A_1}+\ideF{B}$, then $k$ is even.
  \item $\ideF{A_1+B}+\ideF{A_1-B}=\ideF{2\sqrt{p}^{\,t}}\JF^2$ for some $t\in \N_0$ and some ideal  $\JF$ of $\OF$ such that
   $\JF+\ideF{2\sqrt{p}}=\OF$. 
  \item $\IK$ is principal.
 \end{enumerate}
%  Sean $\K=\Q(\sqrt[4]{p})$, $\F=\Q(\sqrt{p})$ con $p\equiv 7\pmod{16}$ un primo racional positivo
%  $A_1=a_1+a_3\sqrt{p}$, $A_2=a_2+a_4\sqrt{p}$
%  y $\alpha=A_1+A_2\sqrt[4]{p}\in\OK-\OF$ tal que
%  $\ideK{\alpha}=\IK^2$, $\nkf{\alpha}=B^2$ y $\alpha$ cumple las condiciones de los casos {\rm 1}, {\rm 2} \'o 
%  {\rm 3} descritos al principio de esta secci\'on. 
%  Entonces se cumplen las siguientes afirmaciones:
%  \begin{enumerate}
%   \item Si $B=b_1+b_2\sqrt{p}$, entonces $b_1,b_2$ tienen la misma paridad en los casos {\rm 1} y {\rm 2}, mientras que en el
%    caso {\rm 3}, $b_1$ es impar y $b_2\equiv 0\pmod{4}$.
%   \item $\ideF{L_2}^2\mid \mid \ideF{A_1+B}+\ideF{A_1-B}$.
%   \item En los casos {\rm 1} y {\rm 2}, $\ideF{A_1+B}+\ideF{A_1-B}=\ideF{L_2}(\ideF{A_1}+\ideF{B})$. En el caso {\rm 3},
%     $\ideF{A_1+B}+\ideF{A_1-B}=\ideF{2}(\ideF{A_1}+\ideF{B})$.
%   \item Sea $\p_\F$  un ideal de $\OF$ tal que $\p_\F\mid \ideF{A_1}$ y $\p_\F\mid \ideF{B}$. Si $\p_\F$ es inerte en $\K/\F$ y $\p_\F^k\mid\mid \ideF{A_1}+\ideF{B}$,
%    entonces $k$ es par.
%   \item Sea $\p_\F$  un ideal de $\OF$ tal que $\p_\F\mid \ideF{A_1}$ y $\p_\F\mid \ideF{B}$. Si $\p_\F$ se descompone 
%    en $\K/\F$ y $\p_\F^k\mid\mid \ideF{A_1}+\ideF{B}$, entonces $k$ es par.
%   \item $\ideF{A_1+B}+\ideF{A_1-B}=\ideF{2\sqrt{p}^{\,t}}\JF^2$ para alg\'un $t\in \N_0$ y $\JF$ un ideal de $\OF$ tal que
%    $\JF+\ideF{2\sqrt{p}}=\OF$. 
%   \item $\IK$ es principal.
%  \end{enumerate}
\end{prop}

 \begin{dem}
  Fitst observe that $$B^2=(a_1+a_3\sqrt{p})^2-\sqrt{p}(a_2+a_4\sqrt{p})^2\quad\mbox{ and }\quad
  B^2=b_1^2+p\, b_2^2+2\,b_1\,b_2\sqrt{p},$$ so:
%  Primero observemos que $$B^2=(a_1+a_3\sqrt{p})^2-\sqrt{p}(a_2+a_4\sqrt{p})^2\quad\mbox{ y }\quad
%  B^2=b_1^2+p\, b_2^2+2\,b_1\,b_2\sqrt{p},$$ as\'i:
 \begin{equation}\label{eqdosformasdeB}
  (a_1^2+a_3^2\,p-2\,p\,a_2\,a_4)+\sqrt{p}(2\,a_1\,a_3-a_2^2-p\,a_4^2)=b_1^2+p\, b_2^2+2\,b_1\,b_2\sqrt{p}.
 \end{equation}
  If 1 or 2 of Theorem \ref{resultofazpr2} holds, since $a_1,a_3$ are odd and $a_2,a_4$ are even, we have
  $$(a_1^2+a_3^2\,p-2\,p\,a_2\,a_4)=b_1^2+p\, b_2^2\equiv 0\pmod{4},$$
  then, $b_1,b_2$ must have the same parity. So, assertion 1 holds for these two cases, 
  but more can be said about $b_1$ and $b_2$. 
  First, note that $\nk{\alpha}=\nf{B^2}\equiv 4\pmod 8$, hence $\nf{B}\equiv 2\pmod{4}$, therefore
  $b_1,b_2$ are odd. If {\rm 1} of Theorem \ref{resultofazpr2} holds, then $a_1\,a_3\equiv 1\pmod{4}$, so:
  $$2\,a_1\,a_3-a_2^2-p\,a_4^2\equiv 2-4-0\equiv 6\pmod{8}, $$
  this implies that  $b_1\not\equiv b_2\pmod{4}$. If 2 of Theorem \ref{resultofazpr2} holds, $a_1\,a_3\equiv 3\pmod{4}$ and 
  $$2\,a_1\,a_3-a_2^2-p\,a_4^2\equiv 6-0-4\equiv 2\pmod{8}, $$
  therefore, $b_1\equiv b_2\pmod 4$.

 %  En los primeros dos casos, como $a_1,a_3$ son impares y $a_2,a_4$ son pares, tenemos  
%  $$(a_1^2+a_3^2\,p-2\,p\,a_2\,a_4)=b_1^2+p\, b_2^2\equiv 0\pmod{4},$$
%  as\'i, $b_1,b_2$ deben tener la misma paridad. Con esto demostramos la afirmaci\'on 1 en los primeros dos casos,
%  pero podemos decir m\'as acerca de $b_1$ y $b_2$. 
% Primero notemos que $\nk{\alpha}=\nf{B^2}\equiv 4\pmod 8$, entonces $\nf{B}\equiv 2\pmod{4}$, de donde
%  $b_1,b_2$ son impares. En el caso {\rm 1}, $a_1\,a_3\equiv 1\pmod{4}$, as\'i:
%  $$2\,a_1\,a_3-a_2^2-p\,a_4^2\equiv 2-4-0\equiv 6\pmod{8}, $$
%  lo que implica  $b_1\not\equiv b_2\pmod{4}$. En el caso 2, $a_1\,a_3\equiv 3\pmod{4}$ y 
%  $$2\,a_1\,a_3-a_2^2-p\,a_4^2\equiv 6-0-4\equiv 2\pmod{8}, $$
%  por tanto, $b_1\equiv b_2\pmod 4$. 
% 

 Now consider the case 3 of Theorem \ref{resultofazpr2}. We have
 $a_1$ is odd and $a_3\equiv 0\pmod 4$. Also,
 $$b_1^2+p\, b_2^2=(a_1^2+a_3^2\,p-2\,p\,a_2\,a_4)\equiv 1+0-0\pmod{8},$$
 and since $p\equiv 7\pmod 8$, then $b_1$ is odd and $b_2$ is even, more precisely, $b_2\equiv 0\pmod{4}$.

%  Ahora estudiemos el tercer caso. Aqu\'i tenemos
%  $a_1$ impar y $a_3\equiv 0\pmod 4$. Por otra parte,
%   $$b_1^2+p\, b_2^2=(a_1^2+a_3^2\,p-2\,p\,a_2\,a_4)\equiv 1+0-0\pmod{8},$$
%   y debido a que $p\equiv 7\pmod 8$, entonces $b_1$ es impar y $b_2$ es par, m\'as a\'un, $b_2\equiv 0\pmod{4}$.
% 

 Let $f(x)=x^2-2(a_1+a_3\sqrt{p}) x+B^2\in\OF[x]$. Let us observe that $f(\alpha)=0$ and since $\alpha\in\OK-\OF$, 
 using Proposition \ref{3prop13} there is $C\in \OF$ such that
\begin{equation}\label{eqcuadradosocuadradosp}
4A_1^2-4B^2=4\left( {A_1+B} \right) \left( {A_1-B} \right)=
 C^2\sqrt{p},
\end{equation}
where $A_1=a_1+a_3\sqrt{p}$ and $B=b_1+b_2\sqrt{p}$.

%  Sea $f(x)=x^2-2(a_1+a_3\sqrt{p}) x+B^2\in\OF[x]$. Observemos que $f(\alpha)=0$ y como $\alpha\in\OK-\OF$, 
%  por la Proposici\'on \ref{dicrenOF} existe $C\in \OF$ tal que  
% \begin{equation}\label{eqcuadradosocuadradosp}
% 4A_1^2-4B^2=4\left( {A_1+B} \right) \left( {A_1-B} \right)=
%  C^2\sqrt{p},
% \end{equation}
% donde $A_1=a_1+a_3\sqrt{p}$ y $B=b_1+b_2\sqrt{p}$.

If 1 of Theorem \ref{resultofazpr2} holds, $a_1\equiv a_3\pmod{4}$ and $b_1\not\equiv b_2\pmod{4}$, so $a_1+b_1\not\equiv a_3+b_2\pmod{4}$
and $a_1-b_1\not\equiv a_3-b_2\pmod{4}$. Then, if
$A_1+B=c_1+c_2\sqrt{p}$, we have $c_1\not\equiv c_2\pmod{4}$ and both are even. 
Clearly, $2\mid A_1+B$. On the other hand
$\dfrac{c_1+c_2\sqrt{p}}{2}$ has an odd coefficient and an even coefficient,
hence $\nf{\dfrac{c_1+c_2\sqrt{p}}{2}}$ is odd and
$2\,L_2 \nmid A_1+B$. In the same way $2\mid A_1-B$ but $2\,L_2\nmid A_1-B$. In case 2 of Theorem \ref{resultofazpr2},  
$a_1\not\equiv a_3\pmod{4}$ and $b_1\equiv b_2\pmod{4}$,
so again $a_1+b_1\not\equiv a_3+b_2\pmod{4}$ and assertion 2 holds for this case.
If assertion 3 of Theorem \ref{resultofazpr2} holds , we have two possibilities, if $a_1\equiv b_1\pmod{4}$, then $A_1+B\equiv 2+0\sqrt{p}\pmod{4}$ and
$A_1-B\equiv 0+0\sqrt{p}\pmod{4}$. On the other hand, if $a_1\not \equiv b_1\pmod 4$, then
$A_1+B\equiv 0+0\sqrt{p}\pmod{4}$ and $A_1-B\equiv 2+0\sqrt{p}\pmod{4}$. Hence $L_2^2\mid\mid A_1\pm B$ and $L_2^4\mid A_1 \mp B$
where the signs are chosen depending if $a_1$ and $b_1$ are equal or not modulo $4$. 
Therefore $\ideF{L_2}^2\mid\mid \ideF{A_1+B}+\ideF{A_1-B}$, concluding the proof of assertion 2.

% 
% En el caso 1, $a_1\equiv a_3\pmod{4}$ y $b_1\not\equiv b_2\pmod{4}$, as\'i que $a_1+b_1\not\equiv a_3+b_2\pmod{4}$
% y $a_1-b_1\not\equiv a_3-b_2\pmod{4}$. Consecuentemente, si
% $A_1+B=c_1+c_2\sqrt{p}$, entonces $c_1\not\equiv c_2\pmod{4}$ y ambos son pares. 
% Es claro que $2\mid A_1+B$. Por otra parte
% $\dfrac{c_1+c_2\sqrt{p}}{2}$ tiene un coeficiente par y uno impar, entonces $\nf{\dfrac{c_1+c_2\sqrt{p}}{2}}$ es impar y
% $2\,L_2 \nmid A_1+B$. De la misma forma $2\mid A_1-B$ pero $2\,L_2\nmid A_1-B$. En el caso 2,  
% $a_1\not\equiv a_3\pmod{4}$ y $b_1\equiv b_2\pmod{4}$,
% por lo que, de nuevo $a_1+b_1\not\equiv a_3+b_2\pmod{4}$. Tambi\'en en este caso se cumple la afirmaci\'on 2.
% Para el caso 3,
%   tenemos dos posibilidades, si $a_1\equiv b_1\pmod{4}$, entonces $A_1+B\equiv 2+0\sqrt{p}\pmod{4}$ y
%   $A_1-B\equiv 0+0\sqrt{p}\pmod{4}$. Por otra parte, si $a_1\not \equiv b_1\pmod 4$, entonces
%   $A_1+B\equiv 0+0\sqrt{p}\pmod{4}$ y $A_1-B\equiv 2+0\sqrt{p}\pmod{4}$. As\'i que  $L_2^2\mid\mid A_1\pm B$ y $L_2^4\mid A_1 \mp B$
%   donde los signos son elegidos dependiendo de que $a_1$ y $b_1$ sean iguales o distintos m\'odulo $4$. 
%   Por lo tanto $\ideF{L_2}^2\mid\mid \ideF{A_1+B}+\ideF{A_1-B}$, concluyendo la demostraci\'on de la afirmaci\'on 2.
% 

To prove 3, we must observe that $2\,A_1,2\,B\in \ideF{A_1+B}+\ideF{A_1-B}$ and $A_1+B,A_1-B\in \ideF{A_1}+\ideF{B}$, so
$$2(\ideF{A_1}+\ide{B})\subseteq\ideF{A_1+B}+\ideF{A_1-B}\subseteq \ideF{A_1}+\ide{B}.$$ This shows that
$\ideF{A_1+B}+\ideF{A_1-B}=L_2^r(\ideF{A_1}+\ide{B})$ with $0\leq r\leq 2$. In the first two cases of Theorem \ref{resultofazpr2},
since $\nf{B}\equiv 2\pmod 4$, then $r>0$.
Furthermore, $a_1$ and $a_3$ have the same parity, hence $\nf{A_1}$ is even, so $\ideF{L_2}\mid\mid \ideF{A_1}+\ideF{B}$ and
$r=1$. If assertion 3 of Theorem \ref{resultofazpr2} holds, $\nk{\alpha}=\nf{B^2}$ is odd, therefore $r=2$.

% Para la afirmaci\'on 3, observemos primero que $2\,A_1,2\,B\in \ideF{A_1+B}+\ideF{A_1-B}$ y $A_1+B,A_1-B\in \ideF{A_1}+\ideF{B}$, as\'i 
% $$2(\ideF{A_1}+\ide{B})\subseteq\ideF{A_1+B}+\ideF{A_1-B}\subseteq \ideF{A_1}+\ide{B}.$$ Esto nos muestra que 
% $\ideF{A_1+B}+\ideF{A_1-B}=L_2^r(\ideF{A_1}+\ide{B})$ con $0\leq r\leq 2$. En los casos 1 y 2,
% como $\nf{B}\equiv 2\pmod 4$, entonces $r>0$.
% Adem\'as, $a_1$ y $a_3$ tienen la misma paridad, por lo que $\nf{A_1}$ es par, as\'i $\ideF{L_2}\mid\mid \ideF{A_1}+\ideF{B}$ y
% $r=1$. En el caso 3, $\nk{\alpha}=\nf{B^2}$ es impar y por tanto $r=2$. 
% 

Now, we are going to prove assertion 4. Since $\p_\F$ is inert, then $\p_\K=\ideK{\p_\F}$ is a prime ideal in $\OK$. 
Let us suppose that $\p_\F^k\mid\mid \ideF{A_1}+\ideF{B}$ and
$\p_\K^t\mid\mid \ideK{\alpha}$. Given that $B^2=A_1^2-\sqrt{p}\, A_2^2$, $\p_\F^{2k}\mid B^2$
and $\p_\F^{2k}\mid A_1^2$, then $\p_\F^k\mid A_2$. Using this, we have that 
$\p_\K^k\mid \ideK{A_1+\sqrt[4]{p}\, A_2}=\ideK{\alpha}$ and $k\leq t$.
On the other hand, $\p_\K^{2t}\mid\mid\nk{\alpha}=A_1^2-\sqrt{p}\,A_2^2=B^2$ and $\p_\K^t\mid\mid \ideK{B}$. Since the irreducible polynomial of
$\alpha$ in $\F[x]$ is $f(x)=x^2-2\, A_1\, x+B^2$, then $\p_\K^{2t}\mid\alpha^2-2\, A_1\, \alpha+B^2=0$. Using this with
$\p_\K^{2t}\mid B^2$ and $\p_\K^{2t}\mid \alpha^2$, we have that $\p_\K^{2t}\mid 2\,A_1\,\alpha$. As a consequence of $\p_\K^t\mid\mid \alpha$
and $\ideK{2}+\p_\K=\OK$, we have $\p_\K^t\mid \ideK{A_1}$. Since $\p_\K=\ideK{\p_\F}$, then $t\leq k$, hence $t=k$. From the equality
$\ideK{\alpha}=\IK^2$, we know that $t$ must be even, therefore $k$ is even.

% Ahora vamos a demostrar la afirmaci\'on 4. Como $\p_\F$ es inerte, entonces $\p_\K=\ideK{\p_\F}$ es un ideal primo en $\OK$. 
% Supongamos que $\p_\F^k\mid\mid \ideF{A_1}+\ideF{B}$ y
% $\p_\K^t\mid\mid \ideK{\alpha}$. Como $B^2=A_1^2-\sqrt{p}\, A_2^2$, $\p_\F^{2k}\mid B^2$
% y $\p_\F^{2k}\mid A_1^2$, entonces $\p_\F^k\mid A_2$. Por lo anterior, 
% $\p_\K^k\mid \ideK{A_1+\sqrt[4]{p}\, A_2}=\ideK{\alpha}$ y $k\leq t$.
% Por otra parte, $\p_\K^{2t}\mid\mid\nk{\alpha}=A_1^2-\sqrt{p}\,A_2^2=B^2$ y $\p_\K^t\mid\mid \ideK{B}$. Como el polinomio irreducible de 
% $\alpha$ en $\F[x]$ es $f(x)=x^2-2\, A_1\, x+B^2$, entonces $\p_\K^{2t}\mid\alpha^2-2\, A_1\, \alpha+B^2=0$. Lo anterior, 
% $\p_\K^{2t}\mid B^2$ y $\p_\K^{2t}\mid \alpha^2$, implican $\p_\K^{2t}\mid 2\,A_1\,\alpha$. Puesto que $\p_\K^t\mid\mid \alpha$
% y $\ideK{2}+\p_\K=\OK$, entonces $\p_\K^t\mid \ideK{A_1}$. Como $\p_\K=\ideK{\p_\F}$, entonces $t\leq k$ y as\'i $t=k$. De la igualdad
% $\ideK{\alpha}=\IK^2$, tenemos que $t$ debe ser par y por lo tanto $k$ es par.

To prove 5, let us suppose that  $\p_\F$ splits, then $\ideK{\p_\F}=\q_1\q_2$, 
with $\q_1,\q_2$ prime ideals of $\OK$. Let us assume that $\q_1^{2t}\mid\mid \ideK{\alpha}$ 
and $\q_2^{2r}\mid\mid \ideK{\alpha}$. Then $\ideK{\alpha}=\q_1^{2t}\q_2^{2r}\JK$ for some $\JK$ such that $\ideK{\p_\F}+\JK=\OK$.
Since $\nkf{\q_1}=\nkf{\q_2}=\p_\F$, then
$$\nkf{\ideK{\alpha}}=\p_\F^{2(t+r)}\nkf{\JK}=\ideF{B^2},$$ hence $\p_\F^{2(t+r)}\mid\mid B^2$.
Without lost of generality, suppose that $r>t$. Then $r=t+s$ for some $s\in\N$ and 
$$\q_1^{4t}\mid\mid\ideK{\alpha^2},\quad \q_1^{4t+2s}\mid\mid B^2,\quad \q_2^{4t+4s}\mid\mid\alpha^2, \quad 
\q_2^{4t+2s}\mid\mid B^2.$$ 
Using this in the equality
$\alpha^2-2\,A_1\,\alpha+B^2=0$, we have that $\q_1^{4t}\mid 2\, A_1 \, \alpha$ and $\q_2^{4t+2s}\mid 2\, A_1 \, \alpha$. 
Since $\q_1^{2t}\mid\mid \ideK{\alpha}$ and $\q_2^{2t+2s}\mid\mid\ideK{\alpha}$, then $(\q_1\q_2)^{2t}=\ideK{\p_\F}^{2t}\mid \ideK{A_1}$.
We will show that $\ideK{\p_\F}^{2t}\mid\mid \ideK{A_1}$. Let us suppose that $\q_1^{2t+1}\mid \ideK{A_1}$.
Since $\q_1^{2t+1}\mid B$, then $\q_1^{4t+2}\mid B^2=A_1^2-\sqrt{p}\,A_2^2$, which implies that $\q_1^{2t+1}\mid A_2$. From this
it follows that $\q_1^{2t+1}\mid\ideK{\alpha}$, which is not possible. Therefore $\q_1^{2t}\mid\mid\ideK{A_1}$. On the other hand, since $A_1\in\OF$, 
then for each $\q_1$ that divides $A_1$ there must exist $\q_2$ that divides $A_1$, hence $\q_2^{2t}\mid\mid\ideK{A_1}$, therefore
$\p_\F^{2t}\mid\mid\ideF{A_1}$. Also, $\p_\F^{2t+s}\mid\mid \ideF{B}$, so $\p_\F^{2t}\mid \ideF{B}$. Therefore, 
$\p_\F^{2t}\mid\mid\ideF{A_1}+\ideF{B}$, where $k=2t$ as asserted in 5.

From 4 and 5, the only prime ideals that can appear an odd number of times in the factorization of $\ideF{A_1}+\ideF{B}$ are ramified
ideals. In this case, this ideals are $\ideF{\sqrt{p}}$ and $\ideF{L_2}$. Using the equality from assertion 3, we may say the same about
the ideal $\ideF{A_1+B}+\ideF{A_1-B}$. Using $2$, we know that $\ideF{L_2}^2\mid\mid \ideF{A_1+B}+\ideF{A_1-B}$. 
This proves assertion 6.

Finaly, from 6, 
\begin{equation}\label{eqprimosrelativos}
 \ideF{\dfrac{A_1+B}{2\,\sqrt{p}^{\,k}\,\JF^2}}+\ideF{\dfrac{A_1-B}{2\,\sqrt{p}^{\,k}\,\JF^2}}=\OF.
\end{equation}
If we use equation (\ref{eqcuadradosocuadradosp}) as an ideal equality, then
$$
\begin{array}{rcl}
 \ideF{4\left( {A_1+B} \right) \left( {A_1-B} \right)} & = & 
  \ideF{4\sqrt{p}^k\JF^2}^2\ideF{\dfrac{A_1+B}{2\,\sqrt{p}^{\,k}\,\JF^2}}\ideF{\dfrac{A_1-B}{2\,\sqrt{p}^{\,k}\,\JF^2}}\\
 & = &  \ideF{C}^2\ideF{\sqrt{p}},\\
\end{array}$$
which implies:
$$\ideF{\dfrac{A_1+B}{2\,\sqrt{p}^{\,k}\,\JF^2}}\ideF{\dfrac{A_1-B}{2\,\sqrt{p}^{\,k}\,\JF^2}}=
 \ideF{\dfrac{{C}}{4\sqrt{p}^{\,k}\JF^2}}^2\ideF{\sqrt{p}},$$
 where all the ideals in the previous equality are integral ideals.
Using (\ref{eqprimosrelativos}), the ideals from the left side are relatively prime, so one of them must be a square and
the other one is a square times $\ideF{\sqrt{p}}$. Let us suppose that: 
$$\ideF{\dfrac{A_1\pm B}{2\,\sqrt{p}^{\,k}\,\JF^2}}=\JJ_1^2,\quad\quad \ideF{2\,A_1\pm 2\,B}=\ideF{2^2\,\sqrt{p}^{\,k}}\JF^2\JJ_1^2.$$
In this way, if $k$ is even, then $\ideF{2\,A_1\pm 2\,B}=\JJ_2^2$, where $\JJ_2^2$ is the ideal of the right side of the equality
and if $k$ is odd, then there exists $\JJ_2\subseteq \OF$ such that $\ideF{2\,A_1 \mp 2\,B}=\JJ_2^2$. In both cases,$\JJ_2^2$ is a principal ideal and, 
since the class number of $\Q(\sqrt{p})$ is odd, then $\JJ_2$ must be principal, say $\JJ_2=\ideF{D}$.
If $A=2\,A_1$, then $A\pm 2B=D^2\,U$ for some $U\in{\cal U}_\K$, where we can suppose that $U=\pm 1$ or $U=\pm U_\F$.
If $U=1$, using Proposition \ref{3prop14}, 
$\sqrt{\alpha}\in\OK$. If $U=-1$, we have
$\nkf{-\alpha}=B^2$ and $t_{\K/\F}(-\alpha)=-2\,A_1=-A $, with $-A\mp 2\,B=D^2$, hence $\sqrt{-\alpha}\in\OK$.
If $U=\pm U_\F$, then
$\nkf{\alpha \, U_\F}=(B\,U_\F)^2$ and $t_{\K/\F}(\alpha\, U_\F)=2\,A_1 \,U_\F$, so
$A_1\,U_\F+B\,U_\F=\pm(D\,U_\F)^2$. Now we proceed as in the previous cases.
Therefore, there is $\mu\in{\cal U}_\K$  such that $\sqrt{\alpha\, \mu}\in\OK$ and it generates $\IK$.
\qed\bigskip
  
 \end{dem}

The previous result requires that $\alpha\not\in\OF$. If this is not the case, we can multiply $\alpha$
by $\mu_1^2$ where $\nkf{\mu_1}=1$, so the norm is preserved, $\K(\sqrt{\alpha})=\K(\sqrt{\alpha\mu_1^2})$
and $\alpha\mu_1^2\not\in\OF$. Observe that $p\equiv 7\pmod{16}$ is needed since the description 
of ${\cal U}_\K$ given in \cite{azpr} depends on this property of $p$.

% El resultado anterior nos pide $\alpha\not\in\OF$. Si esta condici\'on no se cumple, podemos multiplicar $\alpha$
% por $\mu_1^2$ donde $\nkf{\mu_1}=1$, de esta forma la norma se mantiene, $\K(\sqrt{\alpha})=\K(\sqrt{\alpha\mu_1^2})$
% y $\alpha\mu_1^2\not\in\OF$. Observemos que la condici\'on $p\equiv 7\pmod{16}$ es indispensable pues la descripci\'on
% de ${\cal U}_\K$ depende de esta cualidad de $p$.

\begin{cor}
 Let $\K=\Q(\sqrt[4]{p})$, $\F=\Q(\sqrt{p})$ with $0<p\equiv 7\pmod{16}$ a rational prime number, 
 $\alpha=a_1+a_2\sqrt[4]{p}+a_3\sqrt{p}+a_4\sqrt[4]{p^3}\in\OK-\OF$  such that
 $\ideK{\alpha}=\IK^2$, $\alpha$ satisfies one of the assertions {\rm 1}, {\rm 2} or {\rm 3} of {\rm Theorem \ref{resultofazpr2}}
 and $\L=\K(\sqrt{\alpha})$ with
 $\L\neq \K$. Then, $\L/\K$ is a ramified extension or $\L=\K(\sqrt{\mu})$ for some $\mu\in{\cal U}_\K$.\qed
\end{cor}

% 
% 
% 
% \begin{cor}
%    Sean $\K=\Q(\sqrt[4]{p})$, $\F=\Q(\sqrt{p})$ con $p\equiv 7\pmod{16}$ un primo racional positivo, 
%  $\alpha=a_1+a_2\sqrt[4]{p}+a_3\sqrt{p}+a_4\sqrt[4]{p^3}\in\OK-\OF$  tal que
%  $\ideK{\alpha}=\IK^2$, $\alpha$ cumple las condiciones de los casos  {\rm 1}, {\rm 2} \'o {\rm 3} descritos 
%  al principio de la secci\'on y $\L=\K(\sqrt{\alpha})$ con
%  $\L\neq \K$. Entonces, $\L/\K$ es una extensi\'on ramificada \'o $\L=\K(\sqrt{\mu})$ para alguna $\mu\in{\cal U}_\K$.\qed
% \end{cor}

Finally, we will prove the main result.

% 
% 
% Finalmente, usando todo lo que hemos hecho en este cap\'itulo tenemos:
% 
% \begin{teo}\label{2grupodeorden2}
%  Sean $\K=\Q(\sqrt[4]{p})$, $\F=\Q(\sqrt{p})$ con $p\equiv 7\pmod {16}$ un primo racional positivo y
%  $\G2\subseteq\GK$ el $2$-grupo de clases de ideales de $\K$. Entonces $\G2\cong \Z/2\Z$.
% \end{teo}
% 

% \begin{dem}

\bigskip

\noindent {\bf Proof of Theorem \ref{mainresult}.} 
 Let $\L=\K(\sqrt{\alpha})$ for some $\alpha\in\K$ such that $\L\neq \K$. The previous corollary shows that if $\alpha$ is not a unit, then 
 $\L/\K$ is a ramified extension.
 If $\alpha$ is a unit, all the extensions  $\L/\K$ are ramified except
 $\K(\sqrt{U_\F})$. Hence the
 $2$-rank of $\GK$ is 1.
 
 Now we will prove that the order of the $2$-class group is $2$.
 Let $\p_2$ be the only ideal of $\OK$ with $\nk{\p_2}=2$ found in Proposition \ref{3prop11}. 
 Using Proposition \ref{3prop12} we know that $\p_2$ is non-principal but $\p_2^2=\ideK{L_2}$, so $\overline{\p_2}$ is the only class of order $2$ of 
 $\G2$. Let us suppose that there is an ideal  $\IK\subseteq\OK$ such that $\overline{\IK}^2=\overline{\p_2}$.
 Since $\overline{\p_2}$ is it's own inverse, we have $\overline{\IK}^2\,\overline{\p_2}=\overline{\OF}$,
 hence $\II_\K^2\,\p_2$ is a principal ideal.
 We can suppose that $\nk{\IK}$ is odd, because if it is even, $\p_2\mid \IK$, which implies, $\IK=\p_2\,\JK$. Then, 
 $\II_\K^2\,\p_2=\p_2^2\,\JJ_\K^2\,\p_2=\ideK{L_2}\JJ_\K^2\,\p_2,$ so $\JJ_\K^2\p_2$ is related with
 $\II_\K^2\p_2$ and $\nk{\JK}=\dfrac{\nk{\IK}}{2}$. 
 If $\nk{\IK}$ is odd, then $\nk{\IK}^2\equiv 1,9\pmod{16}$, so $\nk{\IK^2\,\p_2}\equiv 2\pmod{16}$.
 As a consequence of this, there must be an element in $\OK$ with norm $\pm 2$, which is not possible by Propositions \ref{3prop11} and \ref{3prop12}.
 Then, there is no $\IK$ such that $\overline{\IK}^2=\overline{\p_2}$ and therefore $\G2\cong \Z/2\Z$. 
 \qed\bigskip

In the next table we give the first rational positive prime numbers $p\equiv 7\pmod{16}$ and the class number of
$\K=\Q(\sqrt[4]{p})$. This values where obtained using the software SAGE  \cite{sage}.

\medskip

% En la siguiente tabla damos los primeros  primos racionales positivos $p\equiv 7\pmod{16}$ 
% y el valor correspondiente de $h_\K$.
% 

$$\begin{tabular}{||c|c||c|c||c|c||c|c||}
 \hline\hline
  $p$ & $h_\K$ & $p$ &$h_\K$ & $p$ & $h_\K$ & $p$ & $h_\K$ \\
 \hline\hline
 7  & 2 & 503  & 2   & 1063  & 2  & 1831  & 6 \\
 \hline
 23 & 2  & 599  & 2  & 1223  & 42  & 1847  & 6 \\
 \hline
 71 & 2  & 631  & 2  & 1303  & 6  & 1879  & 6 \\
 \hline
 103 & 2  & 647  & 2  & 1319  & 2  & 2039  & 2 \\
 \hline
 151 & 2  & 727  & 330  & 1367  & 6  & 2087  & 2 \\
 \hline
 167 & 2  & 743  & 2  & 1399  & 2  & 2311  & 2 \\
 \hline
 199 & 2  & 823  & 2  & 1447  & 2  & 2423  & 6 \\
 \hline
 263 & 2  & 839  & 18  & 1511  & 2  & 2503  & 2 \\
 \hline
 311 & 2  & 919  & 2  & 1543  & 154  & 2551  & 2 \\
 \hline
 359 & 6  & 967  & 2  & 1559  & 2  & 2647  & 2 \\
 \hline
 439 & 50  & 983  & 2  & 1607  & 6  & 2663  & 2 \\
 \hline
 487 & 2  & 1031  & 2  & 1783  & 2  & 2711  & 6 \\
 \hline\hline
\end{tabular}$$

\bigskip

\begin{cor}
 Let $\K=\Q(\sqrt[4]{p})$ with $0<p\equiv 7\pmod{16}$ a rational prime and $h_\K=2$. Then $\H_\K=\K(\sqrt{2})$.
\end{cor}
\begin{dem}
 The only non-ramified quadratic extension of $\K$ is $\K(\sqrt{U_\F})$, then this is the Hilbert class field. The 
 assertion follows from the equality $2=U_\F L_2^2$.\qed\bigskip
\end{dem}

% $$\begin{tabular}{||c|c||c|c||c|c||c|c||}
%  \hline\hline
%  \quad $p$ \quad \quad &\quad  $h_\K$ \quad \quad &\quad  $p$ \quad \quad &\quad  $h_\K$ \quad \quad & \quad $p$ \quad \quad & \quad $h_\K$ \quad \quad & \quad $p$ \quad \quad & \quad $h_\K$ \quad \quad \\
%  \hline\hline
%  7  & 2 & 503  & 2   & 1063  & 2  & 1831  & 6 \\
%  \hline
%  23 & 2  & 599  & 2  & 1223  & 42  & 1847  & 6 \\
%  \hline
%  71 & 2  & 631  & 2  & 1303  & 6  & 1879  & 6 \\
%  \hline
%  103 & 2  & 647  & 2  & 1319  & 2  & 2039  & 2 \\
%  \hline
%  151 & 2  & 727  & 330  & 1367  & 6  & 2087  & 2 \\
%  \hline
%  167 & 2  & 743  & 2  & 1399  & 2  & 2311  & 2 \\
%  \hline
%  199 & 2  & 823  & 2  & 1447  & 2  & 2423  & 6 \\
%  \hline
%  263 & 2  & 839  & 18  & 1511  & 2  & 2503  & 2 \\
%  \hline
%  311 & 2  & 919  & 2  & 1543  & 154  & 2551  & 2 \\
%  \hline
%  359 & 6  & 967  & 2  & 1559  & 2  & 2647  & 2 \\
%  \hline
%  439 & 50  & 983  & 2  & 1607  & 6  & 2663  & 2 \\
%  \hline
%  487 & 2  & 1031  & 2  & 1783  & 2  & 2711  & 6 \\
%  \hline\hline
% \end{tabular}$$
% 

\section{Principal and non-principal ideals}

In this section we will give a criterion to decide if an ideal of $\OK$ is principal or not, in the spirit of Theoremn 18 of \cite{azpr3} and Theorem 
4.4 of \cite{azpr4}.

\begin{teo}\label{principalynoprincipal}
 Let $\K=\Q(\sqrt[4]{p})$ with $0<p\equiv 7\pmod{16}$ a rational prime number and $\IK$ an ideal of $\OK$ such that
 $\mcd{\nk{\IK},2}=1$. The order of the class $\overline{\IK}$ in $\GK$ is odd if and only if 
 $\nk{\IK}\equiv \pm 1\pmod{8}$.
\end{teo}
\begin{dem} We are going to construct an ideal $\IK$ with $\nk{\IK}\equiv 3\pmod{8}$ and we will use it to prove
the assertion. Let $q$ be a rational prime such that $q\equiv 3\pmod{8}$ and $q\equiv a\pmod{p}$, where 
$a\in\Z$ is such that $\leg{a}{p}=-1$. We can guarantee
the existence of such a prime number using Dirichlet's Theorem on infinite primes in an arithmetic sequence,
the Chinese Reminder Theorem and the fact that half the numbers between $1$ and $p-1$ are non-quadratic residues modulo
$p$.

Since $p,q\equiv 3\pmod{4}$ then, $\leg{p}{q}=-\leg{q}{p}$ and the fact that $q\equiv a\pmod{q}$ implies that
$\leg{q}{p}=\leg{a}{p}=-1$. Hence, there exists $b\in \Z$ such that $b^2\equiv p\pmod{q}$. Since $\leg{-1}{q}=-1$
then $\leg{b}{q}=1$ or $\leg{-b}{q}=1$. In both cases, there is $c\in \Z$ such that $c^4\equiv p\pmod{q}$.

As a consequence of the fatorization $x^4-p\equiv (x-c)(x+c)(x^2+c^2)\pmod{q}$ and Dedekind's Theorem 
on the factorization of primes in monogenic 
number fields (\cite{alacawiliams}, Theorem 10.3.1), there are at least two ideals with norm $q$: 
$\II_1=\ideK{q,\sqrt[4]{p}-c}$, $\II_2=\ideK{q,\sqrt[4]{p}+c}$. Let $\alpha\in\OK$, since 
 $$\begin{array}{rcl}
  \nk{\alpha} & = & a_1^4 - a_2^4\, p + 4\, a_1\, a_2^2\, a_3\, p - 2\, a_1^2\, a_3^2\, p-4\, a_1^2\, a_2\, a_4\, 
  p +a_3^4\, p^2\\
  &   & - 4\, a_2\, a_3^2\, a_4\, p^2 + 2\, a_2^2\, a_4^2\, p^2 + 4\, a_1\, a_3\, a_4^2\, p^2 - a_4^4\, p^3,\\\\
 \end{array}$$
 then $\nk{\alpha}$ must be a quartic power modulo $p$, this is $\nk{\alpha}\equiv \pm 1\pmod{8}$.  Hence 
 $\II_1$ is a non-principal ideal. If the order of $\overline{\II_1}\in{\GK}$ is odd, say $k$, then
 $\II_1^k$ is also a non-principal ideal, since $\nk{\II_1^k}\equiv 3\pmod{8}$, therefore,
 the order of $\overline{\II_1}\in{\GK}$ is odd and $\nk{\II_1}\equiv 3\pmod{8}$.
 
 Take $\JK\subseteq\OK$ an ideal in $\OK$ with $\mcd{\nk{\JK},2}=1$. If $\nk{\JK}\equiv \pm 3\pmod{8}$,
 clearly, the order of $\overline{\JK}\in\GK$ is even. Now, suppose that 
 the order of $\overline{\JK}$ is even and $\nk{\JK}\equiv \pm 1\pmod{8}$. Since $h_\K\equiv 2\pmod{4}$, then
 the order of $\overline{\JK}\,\overline{\II_1}$ is odd and $\nk{\JK\II_1}\equiv \pm 3\pmod{8}$, a contradiction.
 Therefore, the order of $\overline{\JK}\in\GK$ is odd.\qed\bigskip
\end{dem}

\begin{cor}
 Let $\K=\Q(\sqrt[4]{p})$ with $0<p\equiv 7\pmod{16}$ a rational prime number and $h_\K=2$, then 
 an ideal $\IK\subseteq \OK$ such that $\mcd{\nk{\IK},2}=1$ is principal if and only if $\nk{\IK}\equiv \pm{1}\pmod{8}$.
 \qed
\end{cor}

\noindent{\bf Example.} Let $\K=\Q(\sqrt[4]{7})$, a number field with class number $2$, hence, 
$\H_\K=\Q(\sqrt[4]{7},\sqrt{2})$. The ideal $\ideK{3}$ factors as:
$$\ideK{3}=\ideK{2+\sqrt{7}}\ideK{3,1+\sqrt[4]{7}}\ideK{3,1-\sqrt[4]{7}}.$$
The ideal $\ideK{2+\sqrt{7}}$ is principal since it has norm $$\nk{\ideK{2+\sqrt{7}}}=9\equiv 1\pmod{8}.$$
The other two ideals are non-principal and 
$$\nk{\ideK{3,1+\sqrt[4]{7}}}=\nk{\ideK{3,1-\sqrt[4]{7}}}=3.$$

If $a$ is an odd rational integer, then $\nk{\ideK{a}}=a^4\equiv 1\pmod{8}$, this is in accordance with
Theorem \ref{principalynoprincipal}.

\end{document}